\documentclass[twoside,a4paper,12pt,leqno]{amsart}

\usepackage{verbatim,amsxtra,amssymb,amsmath,color,psfrag,graphicx}

\definecolor{mygrey}{gray}{0.70}
\definecolor{mygreen}{rgb}{0,.75,0}
\definecolor{myred}{rgb}{1,0,0}
\definecolor{orange}{rgb}{1,.5,0}

\numberwithin{equation}{section}

\newtheorem{thm}[equation]{Theorem}
\newtheorem{conj}[equation]{Conjecture}

\newtheorem{prop}[equation]{Proposition}

\theoremstyle{definition}
\newtheorem{defn}[equation]{Definition}

\theoremstyle{remark}
\newtheorem{rem}[equation]{Remark}

\newcommand{\thmref}[1]{Theorem~\ref{#1}}
\newcommand{\propref}[1]{Proposition~\ref{#1}}

\newcommand{\defref}[1]{Definition~\ref{#1}}
\newcommand{\remref}[1]{Remark~\ref{#1}}

\newcommand{\secref}[1]{Section~\ref{#1}}
\newcommand{\conjref}[1]{Conjecture~\ref{#1}}

\newlength\iii
\newcommand{\num}[1]{\makebox[\iii][r]{(#1)}}

\newcommand\1{\boldsymbol 1}

\newcommand\ab{\boldsymbol a}

\newcommand\be{\beta}
\newcommand\bnd{{\mathop{\text{\bf bnd}\,}}}

\renewcommand\d{\delta}

\newcommand\e{\boldsymbol{\varepsilon}}

\newcommand\f{\varphi}

\newcommand{\G}{\mathbf{G}}
\newcommand{\Ga}{\Gamma}
\newcommand{\ga}{\gamma}
\newcommand\Go{{\mathbf{G}}^{(0)}}
\newcommand\Gd{{\mathbf{G}}^{(2)}}
\newcommand\g{\mathbf{g}}

\newcommand\h{\mathbf{h}}

\renewcommand\l{\lambda}
\newcommand\lsd{\ltimes}

\newcommand\Mor{\operatorname{\mathsf{Mor}}}
\newcommand\Mr{\mathcal{M}}

\newcommand\Obj{\operatorname{\mathsf{Obj}}}

\newcommand\om{\omega}
\newcommand\ov{\overline}

\renewcommand\P{\mathbf{P}}
\renewcommand\part{\partial}
\newcommand\PPi{\boldsymbol \Pi}
\renewcommand\Pr{\mathcal{P}}

\newcommand\R{\mathbb{R}}

\newcommand\s{\operatorname{\mathsf{s}}}
\newcommand\si{\sigma}
\newcommand\Stab{\operatorname{Stab}}

\renewcommand\t{\operatorname{\mathsf{t}}}
\newcommand\toitself{\hookleftarrow}
\newcommand{\toto}{\mathop{\,\longrightarrow\,}}
\renewcommand\th{\theta}

\newcommand\wh{\widehat}
\newcommand\wt{\widetilde}

\newcommand\Z{\mathbb{Z}}

\begin{document}

\title{Amenability and the Liouville property}

\dedicatory{Dedicated to Hillel Furstenberg}

\author{Vadim A. Kaimanovich}

\address{CNRS UMR 6625, IRMAR, Universit\'{e} Rennes-1, Campus Beaulieu, 35042
Rennes, France}

\email{kaimanov@univ-rennes1.fr}

\thanks{}



\begin{abstract}
We present a new approach to the amenability of groupoids (both in the measure
theoretical and the topological setups) based on using Markov operators. We
introduce the notion of an invariant Markov operator on a groupoid and show
that the Liouville property (absence of non-trivial bounded harmonic functions)
for such an operator implies amenability of the groupoid. Moreover, the
groupoid action on the Poisson boundary of any invariant operator is always
amenable. This approach subsumes as particular cases numerous earlier results
on amenability for groups, actions, equivalence relations and foliations. For
instance, we establish in a unified way topological amenability of the boundary
action for isometry groups of Gromov hyperbolic spaces, Riemannian symmetric
spaces and affine buildings.
\end{abstract}

\maketitle

\setcounter{tocdepth}{1}

{\small \tableofcontents}

\thispagestyle{empty}

\section*{Introduction}

The notion of \emph{amenability} for groups is, from the analytical point of
view, the most natural generalization of finiteness or compactness.
\emph{Amenable groups} are those which admit an \emph{invariant mean} (rather
than an invariant probability measure, which is the case for compact groups).
There are many other equivalent definitions of an amenable group, see
\cite{Greenleaf69}, \cite{Paterson88}, \cite{Pier84}. Among the most
constructive is the one formulated in terms of existence of
\emph{approximatively invariant sequences of probability measures} on the group
(\emph{Reiter's condition}), whereas one of the main applications of
amenability is the \emph{fixed point property} for affine actions of amenable
groups on compact spaces.

It turns out that non-amenable groups may still have actions which look like
actions of amenable groups. This observation led Zimmer \cite{Zimmer77},
\cite{Zimmer78} to introduce the notion of an \emph{amenable action}. In the
same way as with groups, there are several definitions of an amenable action.
In particular, amenable actions can be characterized both in terms of a fixed
point property (this was the original definition of Zimmer) and in terms of
existence of a sequence of approximatively equivariant maps from the action
space to the space of probability measures on the group (this is an analogue of
Reiter's condition). Yet another generalization is the notion of amenability
for \emph{equivalence relations} and \emph{foliations}
\cite{Connes-Feldman-Weiss81}. All these objects can be considered as
\emph{measured groupoids}, and the notion of amenability in each particular
case is a specialization of the general notion of an \emph{amenable measured
groupoid}. There is also a similar notion of an \emph{amenable topological
groupoid} (defined in the topological rather than measure theoretical category)
as well. The general references for the theory of amenable groupoids are
\cite{Renault80}, \cite{Anantharaman-Renault00},
\cite{Corlette-Lamoneda-Iozzi02}.

\medskip

A \emph{groupoid} $\G$ is a small category in which each morphism is an
isomorphism, so that it is determined by its \emph{set of morphisms} (also
denoted $\G$) fibered over the \emph{set of objects} $\Go$ via the source $\s$
and the target $\t$ maps. If $\G$ and $\Go$ are topological spaces and the
structure maps are continuous, then $\G$ is called a \emph{topological
groupoid}. Similarly, if $\G,\Go$ are Borel spaces and the structure maps are
Borel, then $\G$ is a \emph{Borel groupoid}. The basic examples of groupoids
are those associated with groups (of course!), group actions and equivalence
relations.

The fibers $\G_x,\G^x$ of the source and the target maps, respectively, are
moved around by morphisms; for instance, $\g\G^{\s(\g)}=\G^{\t(\g)}$ for any
$\g\in\G$. A \emph{Haar system} (defined by analogy with Haar measures on
groups) on a Borel groupoid $\G$ is a Borel system $\l=\{\l^x\}$ of measures on
the fibers of the target map which is \emph{invariant} in the sense that
$\g\l^{\s(\g)}=\l^{\t(\g)}\;\forall\,\g\in\G$. In order to define a
\emph{measured groupoid} in addition to a Haar system $\l$ one has to specify a
Borel measure $\mu$ on the space of objects $\Go$ which is
\emph{quasi-invariant} with respect to the system $\l$. The latter means that
the \emph{global measure} $\l\star\mu$ on $\G$ obtained by integrating the
fiberwise measures $\l^x$ against the measure $\mu$ on the base is
quasi-invariant with respect to the flip $\g\mapsto\g^{-1}$.

Groupoids with a finite Haar system are similar to compact groups
\cite{Hahn78}. In the same way as in the group case, the immediate
generalization of the existence of a finite Haar system is existence of an
\emph{approximatively invariant} sequence $\th_n=\{\th_n^x\}$ of probability
measures on the fibers $\G^x$, i.e., such that
$\bigl\|\g\th_n^{\s(\g)}-\th_n^{\t(\g)}\bigr\|\to 0$. Groupoids with this
property are called \emph{amenable}. Of course, some further assumptions have
to be made and the notion of convergence itself has to be specified depending
on whether we deal with measured or topological groupoids. For measured
groupoids one requires the systems $\th_n$ to be Borel and absolutely
continuous with respect to the Haar system $\l$, and the convergence to be
weak$^*$ in the space $L^\infty(\G,\l\star\mu)$, whereas for topological
groupoids one requires the systems $\th_n$ to be continuous, and the
convergence to be uniform on compact subsets of $\G$.

The definition of an amenable measured groupoid in terms of approximatively
invariant sequences was given by Renault shortly after the work of Zimmer
\cite[Lemma~II.3.4]{Renault80}. However, until recently it remained relatively
unknown to specialists working on amenability of equivalence relations and
group actions in the measure theoretical setup, where either the original
definition of Zimmer or the definition in terms of existence of a
$\G$-invariant mean $L^\infty(\G)\to L^\infty\bigl(\Go\bigr)$ were used, see
\cite{Connes-Feldman-Weiss81}, \cite{Connes-Woods89},
\cite{Adams-Elliott-Giordano94}; in the context of equivalence relations it was
reintroduced in \cite{Kaimanovich97} (at the time I was not aware of Renault's
work). Being very constructive, Renault's definition significantly simplifies
and clarifies proofs of amenability. Compare, for instance, the original proof
of amenability of the boundary action for isometry groups of Gromov hyperbolic
spaces due to Adams \cite{Adams94}, \cite{Adams96} with the recent much shorter
arguments in \cite{Germain00} and \cite{Kaimanovich03}.

Another advantage of the definition of Renault is that it can easily be adapted
to the topological setup, where it found numerous applications to the theory of
$C^*$-algebras, see \cite{Anantharaman87}, \cite{Higson00},
\cite{Higson-Roe00}, \cite{Anantharaman02}, \cite{Valette02},
\cite{Chabert-Echterhoff-Oyono03}.

\medskip

The classical \emph{Liouville theorem} asserts absence of bounded harmonic
functions on the Euclidean space. The notion of a \emph{harmonic function} (one
which satisfies the mean value property with respect to the transition
probabilities) naturally carries over to general Markov operators acting on a
measure space. Such an operator is said to have the \emph{Liouville property}
if it has no non-constant bounded harmonic functions. The link between the
Liouville property and amenability is based on the so-called \emph{\mbox{0--2}
laws} for Markov operators due to Derriennic \cite{Derriennic76} (also see
\cite{Kaimanovich92}), which, in particular, assert that the Liouville property
is equivalent to asymptotic independence of $n$-step transition probabilities
of the initial distribution. This is precisely what is needed for proving
amenability by constructing approximatively invariant sequences of probability
measures. Yet another, less constructive, way of connecting the Liouville
property with amenability consists in using \emph{measure-linear means} on
$\Z_+$. Any such mean applied along the sample paths of a Liouville Markov
chain provides a projection from the space of functions on the state space onto
constants which is invariant with respect to all the symmetries of the
operator, see \cite{Connes-Feldman-Weiss81}, \cite{Lyons-Sullivan84},
\cite{Kaimanovich-Fisher98}.

In order to realize this idea we introduce the notion of an \emph{invariant
Markov operator} on a groupoid $\G$ by analogy with invariant Markov operators
on groups (corresponding to the usual random walks on groups). The transition
probabilities of such an operator satisfy the equivariance condition
$\pi_{\g'\g}=\g'\pi_\g$ for any $\g',\g\in\G$ whenever the composition $\g'\g$
is well-defined. Invariant Markov operators on $\G$ are in one-to-one
correspondence with systems of probability measures on the fibers of the target
map $\t:\G\to\Go$, and the product of two invariant Markov operators
corresponds to the usual convolution operation for such systems of measures
(or, for their densities with respect to a Haar system in the absolutely
continuous case).

The definition of an invariant Markov operator implies, in particular, that any
transition probability $\pi_\g$ is concentrated on the corresponding fiber
$\G^{\t(\g)}$, so that the sample paths of the associated Markov chain are
confined to the fibers of the target map. Therefore, an invariant Markov
operator on $\G$ can be considered as a $\G$-invariant collection of Markov
operators on the fibers of the target map.

As particular cases this notion includes \emph{covering Markov operators} (for
instance, the ones associated with the \emph{Brownian motion on covering
manifolds}) and, more generally, families of their \emph{conditional operators}
with respect to the Poisson boundary (or its equivariant quotients), the Markov
operators associated with \emph{random walks on equivalence relations} and
\emph{leafwise Brownian motion on foliations} as well as all known models of
\emph{randomization} of the usual random walk on a group (random walks in
\emph{random environment}, with \emph{internal degrees of freedom}, with
\emph{random transition probabilities}).

\medskip

We require invariant Markov operators on groupoids to act on the fiberwise
$L^\infty$ spaces with respect to a Haar system, and call such an operator
\emph{fiberwise Liouville} if its restrictions to the fibers of the target map
have the Liouville property. The main result of the paper is that \emph{if a
groupoid carries a fiberwise Liouville invariant Markov operator then it is
amenable}. We prove it both for measured and topological groupoids
(\thmref{thm:liouville measured} and \thmref{thm:top amenab}, respectively).
Particular cases of this result were earlier established for groups
\cite{Azencott70}, \cite{Furstenberg73}, for equivalence relations and
foliations \cite{Connes-Feldman-Weiss81}, for the Brownian motion on covering
manifolds \cite{Lyons-Sullivan84}, for general covering Markov operators
\cite{Kaimanovich95} as well as for various models of randomization of the
usual random walk on a discrete group, see \cite{Sunyach87},
\cite{Kaimanovich-Kifer-Rubshtein02}.

It is plausible that the converse is also true, at least for measured
groupoids, namely, that \emph{any amenable groupoid carries a fiberwise
Liouville invariant Markov operator}. This is known to be the case for groups
\cite{Rosenblatt81}, \cite{Kaimanovich-Vershik83} (it had been previously
conjectured by Furstenberg \cite{Furstenberg73}), for discrete equivalence
relations (in view of the Connes--Feldman--Weiss theorem on the coincidence of
amenability and hyperfiniteness \cite{Connes-Feldman-Weiss81}) and for group
actions (in a somewhat weaker form, though; by \cite{Elliott-Giordano93},
\cite{Adams-Elliott-Giordano94} any amenable measure class preserving action of
a locally compact group $G$ can be realized as the action on the Poisson
boundary of an appropriate $G$-invariant operator on the product of $G$ by a
countable set).

As a specialization of our main result to the case of groupoids associated with
group actions we show (once again both in the measure theoretical and
topological setups, \thmref{thm:amenab actions-measurable} and
\thmref{thm:amenab actions}, respectively) that \emph{if there exists an
equivariant map assigning to the points from the action space $X$ minimal
harmonic functions of a certain $G$-invariant Markov operator on another space
$S$ (endowed with a proper action of the group $G$), then the action of $G$ on
$X$ is amenable}. The reason is that the \emph{Doob transforms} of the original
operator determined by these minimal harmonic functions have the Liouville
property. In the particular case of the random walk on a countable group with a
finitely supported transition probability this result was also independently
obtained in \cite{Biane-Germain02}.

Typically, such a situation arises when $X$ is a certain \emph{boundary} of the
space $S$. Identification of the space of minimal harmonic functions of a
$G$-invariant Markov operator (in other words, of the \emph{Poisson boundary}
in the measure theoretical setup or the \emph{minimal Martin boundary} in the
topological setup) is, in general, a difficult task (e.g., see
\cite{Kaimanovich96}). However, natural geometrical boundaries are known to
produce minimal harmonic functions in several situations of hyperbolic flavour:
for \emph{simply connected negatively curved manifolds with pinched sectional
curvatures} and, more generally, \emph{Gromov hyperbolic spaces}; for
\emph{Riemannian symmetric spaces}; for \emph{affine buildings}. Therefore, our
result immediately implies the strongest possible \emph{topological amenability
of the boundary actions of the groups of isometries} of these spaces
(\thmref{thm:geometric amenability}). In particular, these groups are
\emph{amenable at infinity}. This provides a unified generalization of numerous
earlier results on amenability of boundary actions \cite{Bowen77},
\cite{Vershik78}, \cite{Zimmer84}, \cite{Spatzier87}, \cite{Spatzier-Zimmer91},
\cite{Adams94}, \cite{Adams96}, \cite{Robertson-Steger96},
\cite{Ramagge-Robertson96}, \cite{Cutting-Robertson03}.

Yet another particular case is the \emph{measure theoretical amenability of the
action of a locally compact group $G$ on the Poisson boundary} either of a
usual random walk on $G$ or of a certain $G$-invariant chain on a $G$-space,
which was earlier established in \cite{Zimmer78} and \cite{Connes-Woods89}
(also see \cite{Elliott-Giordano93}, \cite{Adams-Elliott-Giordano94}).

More generally, following the considerations for covering Markov operators in
\cite{Kaimanovich95} we introduce the notion of the \emph{Poisson extension} of
an invariant Markov operator $P$ on a groupoid $\G$. This is the measured
groupoid associated with the action of $\G$ on the Poisson boundary of the
operator $P$. We prove that \emph{the Poisson extension is amenable}
(\thmref{thm:extension}), which subsumes \thmref{thm:liouville measured} and
provides a generalization of the above results on the amenability of the
Poisson boundary in the group case.

\medskip

These results were presented at a number of seminars and conferences
(University of Chicago 1997, ENS Lyon 1998, University of Genova 1999, Rokhlin
memorial conference, St. Petersburg 1999, University of Orleans 2000,
University of Neuch\^{a}tel 2001, Caltech 2002). In particular, at the seminar of
Anantharaman--Delaroche and Renault in Orleans in May 2000 a draft version of
the present article was circulated.

\medskip

The paper has the following structure. In \secref{sec:groupoids} and
\secref{sec:amenability} we recall main definitions concerning groupoids and
their amenability. The exposition here is mostly based on the books
\cite{Renault80} and \cite{Anantharaman-Renault00}. In \secref{sec:operators}
we introduce the notion of an invariant Markov operator on a groupoid and
discuss various examples of such operators. In \secref{sec:liouville} we prove
amenability of measured groupoids with fiberwise Liouville invariant Markov
operators, and in \secref{sec:poisson extension} we prove amenability of the
Poisson extension of an invariant Markov operator on a measured groupoid.
Finally, in \secref{sec:topological} we establish analogues of these results in
the topological category.

\medskip

The spirit of this paper owes a lot to the work and ideas of \emph{Hillel
Furstenberg}. It was him who laid the foundation of the modern probabilistic
boundary theory of algebraic structures, and, in particular, established the
first results relating amenability to the Liouville property for random walks
on groups --- which are the starting point of the present article. I dedicate
this paper to him with admiration.

\section{Groupoids} \label{sec:groupoids}

\subsection{General definitions}

Let us first recall the definition of a \emph{groupoid} as a \emph{small
category in which each morphism is an isomorphism}.

In other words, a groupoid $\G$ is determined by a \emph{set of objects} (also
called the \emph{set of units})
$$
\Obj\G=\Go
$$
and a \emph{set of morphisms}
$$
\Mor\G \cong \G
$$
(which, following the well-established tradition, we shall usually
denote just by $\G$) endowed with the \emph{source} and \emph{target} maps
$$
\s,\t:\G\to\Go \;.
$$
Denote by
$$
\Gd = \{ (g_1,g_2)\in\G\times\G: \s(\g_1)=\t(\g_2) \}
$$
the set of \emph{composable pairs} in $\G$. The composition is a map
$$
\Gd\to\G \;, \qquad (\g_1,\g_2)\mapsto \g_1\g_2
$$
such that
$$
\s(\g_1\g_2) = \s(\g_2) \;,\qquad \t(\g_1\g_2)=\t(\g_1) \;.
$$

\begin{rem}
Our notations match those used for left actions of groups:
$(g_1g_2)x=g_1(g_2x)$, i.e., $g_2$ is applied ``first''. Alternatively, one
could choose the ``postfix'' notation (which also has some euristic advantages)
corresponding to the right actions with $x(g_1 g_2)=(x g_1)g_2$.
\end{rem}

The composition has the usual properties. Namely, there is an embedding
$$
\e:\Go\to\G \;,
$$
which associates to any object $x\in\Go$ the \emph{identical automorphism}
$\e_x$ such that
$$
\s(\e_x) = \t(\e_x) = x \qquad \forall\,x\in\Go
$$
(usually we shall just identify $x$ and $\e_x$). For any $\g\in\G$ there is a
unique \emph{inverse morphism} $\g^{-1}$ with the property that
$$
\s(\g^{-1}) = \t(\g) \;, \qquad \t(\g^{-1}) = \s(\g)
$$
and
$$
\g\g^{-1} = \e_{\t(\g)} \;, \qquad \g^{-1}\g = \e_{\s(\g)} \;.
$$
Finally, the composition (when well-defined) is associative.

The fibers of the source and target maps are denoted
$$
\G_x = \s^{-1}(x) \;, \qquad \G^x = \t^{-1}(x) \;,
$$
respectively. More generally, for any two subsets $X,Y\subset\Go$ and a subset
$A\subset\G$, we put
$$
A_X = A\cap \s^{-1}(X) \;, \qquad A^Y = A\cap \t^{-1}(Y) \;, \qquad A_X^Y =
A_X\cap A^Y \;.
$$

The \emph{isotropy group} of an object $x\in\Go$ is then $\G_x^x$ (whose unit
is $\e_x$; this is the reason why $\Go$ is called the set of units). The set
$$
\G' = \{\g\in\G:\s(\g)=\t(\g) \} = \bigcup_{x\in\Go}\G_x^x
$$
is called the \emph{isotropy bundle}. The groupoid $\G$ determines the
associated \emph{orbit equivalence relation}
$$
R_\G = \{ (\s(\g),\t(\g)): \g\in\G \} \subset \Go\times\Go
$$
on the space of objects $\Go$. Its classes are called \emph{orbits} of $\G$ in
$\Go$. Therefore, any groupoid can be considered as an extension of its orbit
equivalence relation by the isotropy bundle.

\subsection{Examples of groupoids} \label{sec:examples}

The most basic examples are:

\medskip

\num{i} \emph{Groups}. Any group $G$ can be in a trivial way considered as a
groupoid by putting
$$
\G=G \quad \text{and} \quad \Go=\{o\}
$$
for a single point $o$ with $\s(\g),\t(\g)\equiv o$ and the same composition
rule in $\G$ as in the original group $G$. Then clearly $\G$ consists just of
the isotropy group $\G_o^o$ isomorphic to $G$.

\medskip

\num{ii} \emph{Equivalence relations}. If $R\subset X\times X$ is an
equivalence relation on a set $X$, then for the associated groupoid
$$
\G=R \quad \text{and} \quad \Go=X
$$
with source and target maps
$$
\s(x,y)=y \;, \quad \t(x,y)=x \;,
$$
respectively, the composition
$$
(x,y)(y,z)=(x,z) \;,
$$
and the inverse map
$$
(x,y)^{-1}=(y,x) \;.
$$
The embedding $\e:\Go\to\G$ is diagonal, i.e.,
$$
\e_x=(x,x) \;.
$$
All isotropy groups are trivial, whereas the orbit equivalence relation $R_\G$
on the set of objects $\Go\cong X$ is precisely the original equivalence
relation $R$.

\medskip

\num{iii} \emph{Group actions}. Let a group $G$ acts (on the left) on a set
$X$. For the associated groupoid
$$
\G=\{(gx,g,x): g\in G, x\in X\} \quad \text{and} \quad \Go=X
$$
with the source and target maps
$$
\s(gx,g,x)=x \;, \quad \t(gx,g,x)=gx \;,
$$
respectively, the embedding
$$
\e_x=(x,e,x)
$$
(here $e$ is the identity of the group $G$), the inverse map
$$
(gx,g,x)^{-1}=(x,g^{-1},gx) \;,
$$
and the composition
$$
(g_1g_2x,g_1,g_2x)(g_2x,g_2,x)=(g_1g_2x,g_1g_2,x) \;.
$$
The isotropy group $\G_x^x$ is then isomorphic to the stabilizer
$$
\Stab_x=\{g\in G: gx=x\}
$$
of the point $x$ in the group $G$, and the orbit equivalence relation $R_\G$ is
the usual orbit equivalence relation of the action of the group $G$. In
particular, for the trivial action of the group $G$ on the singleton $\{o\}$ we
obtain the groupoid (i) associated with the group $G$. If the action is free,
then we obtain the groupoid (ii) associated with the orbit equivalence relation
of the action of the group $G$ on $X$.

\medskip

See \cite{Renault80}, \cite{Anantharaman-Renault00} and
\cite{Corlette-Lamoneda-Iozzi02} for more examples of groupoids, for instance,
those arising from transformation pseudogroups and foliated manifolds. Some
further examples (together with a description of invariant Markov operators on
them) are also given in \secref{sec:examples of operators}.

\subsection{Homogeneous spaces} \label{sec:homogeneous spaces}

By analogy with the group case one can also define the notion of a $\G$-space
for a groupoid $\G$. A (left) \emph{$\G$-space} consists of a set $X$ endowed
with a \emph{projection map} $\t=\t_X:X\to\Go$ and an \emph{action map}
$$
\{(\g,x)\in\G\times X: \s(\g)=\t(x)\} \to X \;, \qquad (\g,x)\mapsto \g x
$$
with natural properties, more precisely,
$$
\t(\g x) = \t(\g) \;, \quad \e_{\t(x)} x = x \;, \quad \text{and} \quad
\g_1(\g_2 x) = (\g_1\g_2) x
$$
whenever the corresponding products are well-defined.

Any $\G$-space gives rise (in the same way as for usual group actions, see
example (iii) in \secref{sec:examples} above) to the associated groupoid called
the \emph{semi-direct product} of $\G$ and $X$ and denoted $\G\lsd X$. Namely,
$$
\Mor (\G\lsd X) = \{(\g x,\g,x): \s(\g)=\t(x) \} \quad \text{and} \quad \Obj
(\G\lsd X) = X \;,
$$
whereas the structure maps are defined in precisely the same way as in example
(iii) from \secref{sec:examples}. In particular, for any groupoid $\G$ its set
of objects $\Go$ can be in a natural way considered as a $\G$-space, and the
semi-direct product $\G\lsd\Go$ is isomorphic to $\G$.

\subsection{Measured groupoids} \label{sec:measured groupoids}

Any morphism $\g\in\G$ determines the bijection
$$
\G^{\s(\g)}\to\G^{\t(\g)} \;, \qquad \h\mapsto\g\h \;.
$$
This observation prompts one to define (by the analogy with the classical
notion of a Haar measure on a locally compact group) the notion of a \emph{Haar
system of measures} on a groupoid $\G$ as a family of measures $\l=\{\l^x\}$ on
the fibers $\G^x\subset\G$ of the target map $\t:\G\to\Go$ which is
$\G$-\emph{invariant} in the sense that
\begin{equation} \label{eq:Haar}
\g\l^{\s(\g)} = \l^{\t(\g)} \qquad \forall\,\g\in\G \;.
\end{equation}

More precisely, we shall say that a groupoid $\G$ is \emph{Borel} if it is
endowed with a Borel structure such that the structure maps are Borel, where
$\Go$ and $\Gd$ are given the Borel structures induced by $\G$ and
$\G\times\G$, respectively. Then a family of measures $\{\l^x\}$ on $\G$
indexed by points $x\in\Go$ is called \emph{Borel} if for any non-negative
Borel function $f$ on $\G$ the function on $\Go$
\begin{equation} \label{eq:lambda_f}
\l(f): x\mapsto\langle\l^x,f\rangle
\end{equation}
is also Borel. Such a family is called \emph{proper} if there exists a Borel
function $f$ with $\l(f)\equiv 1$. Then a \emph{Borel Haar system} on a Borel
groupoid $\G$ is a $\G$-invariant proper Borel family of measures $\{\l^x\}$
concentrated on the sets $\G^x$ (i.e., $\l^x(\G\setminus\G^x)=0$ for any
$x\in\Go$).

A measure $\nu$ on $\G$ is called \emph{quasi-symmetric} if it is
quasi-invariant with respect to the inverse map $\g\mapsto\g^{-1}$. One can
integrate a Borel system of measures $\l=\{\l^x\}$ on the fibers $\G^x$ of the
target map $\t:\G\to\Go$ with respect to any Borel measure $\mu$ on $\Go$,
which gives the Borel measure $\l\star\mu$ on $\G$. Then the measure $\mu$ is
called \emph{quasi-invariant} with respect to the system $\l$ if the measure
$\l\star\mu$ is quasi-symmetric.

Finally, a \emph{measured groupoid} is a triple $(\G,\l,\mu)$, where $\G$ is a
Borel groupoid, $\l$ is a Borel Haar system, and $\mu$ is a measure on $\Go$
quasi-invariant with respect to $\l$. Actually, in this definition we only need
the class of the measure $\mu$, and below we shall also apply the term
``measured groupoid'' to the situation when just a quasi-invariant measure
class (rather than a specific measure) on the space of objects is given.

Given a Borel groupoid $\G$ endowed with a Borel Haar system $\l$ and a Borel
$\G$-space $X$, the projection map $X\mapsto\Go$ allows one to lift the system
$\l$ to a Borel Haar system (also denoted $\l$) of the semi-direct product
$\G\lsd X$. A Borel measure $\mu$ on $X$ is then called \emph{quasi-invariant}
with respect to $(\G,\l)$ if the measure $\l\star\mu$ on the groupoid $\G\lsd
X$ is quasi-symmetric, i.e., if the triple $(\G\lsd X,\l,\mu)$ is a measured
groupoid \cite[Definition~3.1.1]{Anantharaman-Renault00} (in the case when
$X=\Go$ this definition clearly agrees with the previous definition of a
quasi-invariant measure on the space of objects, see \secref{sec:homogeneous
spaces}).

In particular, any measured groupoid $(\G,\l,\mu)$ naturally acts on itself,
and the measure $\l\star\mu$ on $\G$ is quasi-invariant with respect to
$(\G,\l)$ \cite[Proposition 3.4]{Hahn78} (the associated semi-direct product
$\G\lsd\G$ is isomorphic to $\Gd$).

\subsection{Topological groupoids}

In the \emph{topological setting}, one assumes that the groupoid $\G$ is a
topological space and that the structure maps are continuous, where $\Gd$ has
the topology induced by $\G\times\G$, and $\Go$ has the topology induced by
$\G$. Furthermore, the source and the target maps are surjective and open.
Similar assumptions are made in the definition of a continuous $\G$-space. We
shall be concerned exclusively with topological groupoids and continuous
$\G$-spaces which are \emph{second countable, locally compact and Hausdorff}.

We shall use the standard definition of a \emph{continuous Haar system}
$\lambda$ for a locally compact groupoid $\G$: continuity in this situation
means that for any compactly supported continuous function $f$ on $\G$ the
associated function $\l(f)$ \eqref{eq:lambda_f} on $\Go$ is also continuous.

\section{Amenability} \label{sec:amenability}

\subsection{Amenable groups}

The notion of \emph{amenability} was first introduced for groups and is the
most natural, from the analytical point of view, generalization of finiteness
or compactness. Namely, compact groups are distinguished within the class of
all locally compact topological groups by the property that they carry a
\emph{finite} invariant measure.

Recall that a \emph{mean} on a measure space $(X,m)$ is a positive normalized
(hence, continuous) linear functional on the space $L^\infty(X,m)$, or,
equivalently, a \emph{finitely additive} probability measure on $X$ absolutely
continuous with respect to $m$. In a similar way one defines means with values
in $L^\infty(X',m')$ for any quotient space $(X',m')$ of $(X,m)$.

For an \emph{amenable group} $G$ one requires existence of an \emph{invariant
mean} on $L^\infty(G)$ instead of a probability measure (this is the classical
definition of amenability), or, equivalently, by \emph{Reiter's condition},
existence of a sequence of probability measures $\th_n$ on $G$ which is
\emph{approximatively invariant} in the sense that
$$
\left\| g\th_n - \th_n \right\| \toto_{n\to\infty} 0 \qquad\forall\,g\in G \;.
$$
Reiter's condition is one of the most constructive among a host of other
(equivalent) definitions of amenability, see \cite{Greenleaf69}, \cite{Pier84},
\cite{Paterson88}.

\subsection{Amenable groupoids} \label{sec:amenable groupoids}

The notion of amenability can be naturally generalized to groupoids (including
groups as a particular case). A measured groupoid $(\G,\l,\mu)$ is called
\emph{amenable} if there exists a $\G$-invariant mean
$$
\PPi:L^\infty(\G,\l\star\mu)\to L^\infty(\Go,\mu) \;,
$$
where $\Go$ is considered as the quotient of $\G$ under the target map
\cite[Definition~3.2.8]{Anantharaman-Renault00}. By \emph{$\G$-invariance} we
mean that $\PPi$ commutes with the convolution with any system of finite
measures $\th=\{\th^x\}$ on the fibers of the target map absolutely continuous
with respect to the Haar system, i.e.,
$$
\PPi(\th\ast f) = \th\ast\PPi(f) \qquad \forall\,f\in L^\infty(\G,\l\star\mu)
\;,
$$
where the \emph{convolution} is defined as
$$
\th\ast f(\g) = \int f(\h^{-1}\g) \,d\th^{\t(\g)}(\h) \;.
$$

As in the group case, this definition has a constructive counterpart formulated
in terms of systems of probability measures $\th=\{\th^x\}_{x\in\Go}$ on the
fibers $\G^x\subset\G$ of the target map $\t:\G\to\Go$ which are absolutely
continuous with respect to the Haar system $\l=\{\l^x\}$ and measurable in the
sense that the Radon--Nikodym derivative $d\th^{\t(\g)}/d\l^{\t(\g)}(\g)$ is a
measurable function on $\G$. Amenability of a measured groupoid $(\G,\l,\mu)$
is then equivalent to existence of a sequence $\th_n$ of such systems which is
\emph{approximatively invariant} in the sense that
\begin{equation} \label{eq:def-amen}
\int \left\| \g\th_n^{\s(\g)} - \th_n^{\t(\g)} \right\| f(\g)\,d\l\star\mu(\g)
\toto_{n\to\infty} 0
\end{equation}
for any test function $f\in L^1(\G,\l\star\mu)$ (this is a reformulation of
condition (iv) from \cite[Proposition~3.2.14]{Anantharaman-Renault00}).

\subsection{Topological amenability} \label{sec:top amenab}

In the topological setup, for a locally compact groupoid $\G$, a system of
probability measures $\th=\{\th^x\}$ on the fibers $\G^x\subset\G,\,x\in\Go$ of
the target map $\t:\G\to\Go$ is called \emph{continuous} if for any continuous
function with compact support $f$ the associated function $\th(f)$
\eqref{eq:lambda_f} is continuous on $\Go$. A locally compact groupoid $\G$ is
called \emph{topologically amenable} if it admits a sequence $\th_n$ of
continuous systems of probability measures on the fibers $\G^x$ which is
\emph{topologically approximatively invariant} in the sense that
$$
\left\| \g\th_n^{\s(\g)} - \th_n^{\t(\g)} \right\| \toto_{n\to\infty} 0
$$
uniformly on compacts in $\G$ \cite[Definitions 2.2.1 and
2.2.8]{Anantharaman-Renault00}. In the presence of a continuous Haar system
$\l$ on $\G$ the systems $\th_n$ from the above definition can be chosen in
such a way that the Radon--Nikodym derivatives
$d\th_n^{\t(\g)}/d\l_n^{\t(\g)}(\g)$ are continuous on $\G$ \cite[Proposition
2.2.13]{Anantharaman-Renault00}.

\medskip

See \cite{Anantharaman-Renault00} for more details on the notion of amenability
for measured and topological groupoids.

\section{Invariant Markov operators} \label{sec:operators}

\subsection{Markov chains and operators}

Let $X$ be a Borel space; a family $\pi=\{\pi_x\}$ of Borel probability
measures on $X$ indexed by points $x\in X$ is called Borel if for any
non-negative Borel function $f$ on $X$ the function
$$
\pi(f):x\mapsto\langle\pi_x,f\rangle
$$
is also Borel. Any such family determines a \emph{Markov chain} on $X$ with the
transition probabilities $\pi_x$. As usual, we denote by $(X^{\Z_+},\P_\th)$
the space of \emph{sample paths} $\ov x= (x_0,x_1,\dots)$ of this chain with
the initial ($\si$-finite) distribution $\th$.

It is convenient to talk about a Markov chain in terms of the associated
\emph{Markov operator}
$$
Pf=\pi(f) \;.
$$
Its \emph{dual operator} acts on the space of non-negative Borel measures on
$X$ by the formula
$$
\th P = \int \pi_x\,d\th(x)
$$
(this is a standard notation in the theory of Markov chains, e.g., see
\cite{Revuz84}), so that
$$
\langle \th, Pf \rangle = \langle \th P, f \rangle
$$
for any non-negative Borel function $f$ and measure $\th$ (both sides are
allowed to take infinite values as well). Then the one-dimensional
distributions of the measure $\P_\th$ in the path space are $\th P^n$.

A measure (not necessarily finite!) $\th$ on $X$ is called \emph{invariant} (or
\emph{stationary}) with respect to the operator $P$ if $\th=\th P$ (or,
equivalently, if the associated measure $\P_\th$ in the path space is
shift-invariant), \emph{quasi-invariant} if $\th$ and $\th P$ are equivalent,
and \emph{adapted} if the measure $\th P$ is absolutely continuous with respect
to $\th$ (for the lack of a well-established term; more legitimate if
cumbersome candidates would be ``sub-quasi-invariant'', ``quasi-excessive'' or
``null-preserved''). In the latter case the operator $P$ also acts on the space
$L^\infty(X,\th)$. Alternatively, one can define a Markov operator directly as
an operator on $L^\infty(X,\th)$ \cite{Foguel69}.

\subsection{Invariant Markov operators}

\begin{defn} \label{def:invariant operator}
A Markov operator on a Borel groupoid $\G$ is called \emph{invariant} if its
transition probabilities $\{\pi_\g\}$ satisfy the relation
\begin{equation} \label{eq:equivariance}
\pi_{\g'\g} = \g'\pi_\g \qquad \forall\,(\g',\g)\in\Gd \;.
\end{equation}
\end{defn}

In other words, any transition probability $\pi_\g$ of an invariant Markov
operator is concentrated on the corresponding fiber $\G^{\t(\g)}$, i.e.,
\begin{equation} \label{eq:target preserved}
\pi_\g(\G^{\t(\g)})=1 \;,
\end{equation}
and the relation \eqref{eq:equivariance} holds for any $\g'\in\G_{\t(g)}$.
Thus, an invariant Markov operator $P$ on $\G$ can be considered as a
$\G$-invariant collection of Markov operators $P_x$ on the fibers
$\G^x,\,x\in\Go$.

\begin{prop} \label{prop:extension}
Any Borel system of probability measures on the fibers $\G^x,\,x\in\Go$ of the
target map $\t:\G\to\Go$ uniquely extends to the system of transition
probabilities of an invariant Markov operator on $\G$.
\end{prop}

\begin{proof}
We shall consider a system of probability measures $\pi_x,\,x\in\Go$ on the
fibers $\G^x$ as the transition probabilities from the associated points
$\e_x\in\G$. If we define
\begin{equation} \label{eq:defsystem}
\pi_\g = \g\pi_{\s(\g)} \qquad \forall\,\g\in\G \;,
\end{equation}
then
$$
\pi_{\g'\g} = (\g'\g)\pi_{\s(\g'\g)} = \g' \bigl( \g\pi_{\s(\g)} \bigr) =
\g'\pi_\g \qquad\forall\, (\g',\g)\in\Gd \;,
$$
so that the system \eqref{eq:defsystem} is $\G$-invariant.
\end{proof}

In the measure theoretical setup, when talking about an invariant Markov
operator $P$ on a measured groupoid $(\G,\l,\mu)$ \emph{we shall always assume
that the measure $\l\star\mu$ is $P$-adapted}, which is equivalent to the
measures $\l^x$ from the Haar system being adapted with respect to the
corresponding fiberwise Markov operators $P_x$ for $\mu$-a.e. $x\in\Go$.

Note that we do not impose on invariant Markov operators $P$ any conditions
related to existence of a \emph{$P$-invariant} measure.

\subsection{Absolutely continuous transition probabilities}

Given a Borel groupoid $\G$ with a Borel Haar system $\{\l^x\}$, we shall say
that the transition probabilities $\{\pi_\g\}$ of an invariant Markov operator
$P$ are \emph{absolutely continuous} if
$$
\pi_\g\prec\l^{\t(\g)} \qquad\forall\,\g\in\G
$$
(we use the symbol $\prec$ to denote the absolute continuity of one measure
with respect to another one). In view of \defref{def:invariant operator} and
\propref{prop:extension} this condition is equivalent to the absolute
continuity just of the measures $\pi_x,\,x\in\Go$ with respect to the
corresponding measures~$\l^x$. It is easy to see that the absolute continuity
of the transition probabilities of an invariant Markov operator $P$ on a
measured groupoid $(\G,\l,\mu)$ implies that the measure $\l\star\mu$ is
$P$-adapted. [Actually, it is sufficient to require absolute continuity of
transition probabilities $\pi_\g$ just for $\l\star\mu$-a.e. $\g\in\G$, or,
equivalently, just for $\mu$-a.e. $\e_x\cong x\in\Go$.] Therefore, an invariant
Markov operator with absolutely continuous transition probabilities acts on the
space $L^\infty(\G,\l\star\mu)$.

\begin{rem}
Absolute continuity of transition probabilities is not necessary for the
measure $\l\star\mu$ being adapted with respect to an invariant Markov operator
$P$. The simplest example is provided by random walks on groups, see Example
(i) below, in which case the Haar measure is adapted (even invariant) with
respect to any random walk on the group.
\end{rem}

\begin{rem}
As we have seen in \propref{prop:extension}, invariant Markov operators on a
groupoid $\G$ are in one-to-one correspondence with systems of probability
measures on the fibers of the target map $\t:\G\to\Go$. The product of
invariant Markov operators corresponds then to the usual convolution in the
space of such systems, or, for operators with absolutely continuous transition
probabilities, in the space of densities of these measures with respect to a
Haar system, see \cite{Renault80}, \cite{Anantharaman-Renault00} for the
definition of convolution on groupoids.
\end{rem}

\subsection{Examples of invariant Markov operators} \label{sec:examples of
operators}

Markov operators satisfying reasonable ``homogeneity conditions'' can, as a
rule, be interpreted as invariant Markov operators on appropriate groupoids.

\medskip

\num{i} \emph{Random walks on groups.} According to the classical definition, a
(right) \emph{random walk} on a locally compact group $G$ determined by a
probability measure $\pi$ on $G$ is the Markov chain on $G$ with the transition
probabilities $\pi_g=g\pi$ which are equivariant with respect to the action of
the group on itself on the left. Its transition operator $P$ is an invariant
Markov operator on the associated groupoid $\G$ (see example (i) in
\secref{sec:examples}). For introducing a structure of a measured groupoid on
$\G$ it is sufficient to take a left Haar measure on the group $G$ (which can
be considered as a ``Haar system'', since the object space $\Go$ is a
singleton), and absolute continuity of transition probabilities of the operator
$P$ is equivalent just to absolute continuity of the measure $\pi$ with respect
to $\l$.

\medskip

\num{ii} \emph{Markov chains on equivalence relations and foliations.} An
equivalence relation $R$ on a standard Borel space $X$ is called
\emph{standard} if it is a Borel subset of $X\times X$, and it is called
\emph{countable} if the equivalence class (the \emph{leaf}) $[x]=R(x)=\{y:
(x,y)\in R\}$ of any point $x\in X$ is at most countable. A standard countable
equivalence relation is also called \emph{discrete}. A standard equivalence
relation $R$ is called \emph{non-singular} with respect to a Borel probability
measure $\mu$ on $X$ (or, equivalently, the measure $\mu$ is
\emph{quasi-invariant} with respect to $R$) if for any subset $A\subset X$ with
$\mu(A)=0$ its \emph{saturation} $[A]=\bigcup_{x\in A}[x]$ also has measure 0,
see \cite{Feldman-Moore77}.

Let $\G$ be the groupoid associated with the equivalence relation $R$, see
\secref{sec:examples}. Then the counting measures on the sets $\G^x=x\times
[x]$ provide a Haar system $\l$ for $\G$, and the $R$-quasi-invariance of a
measure $\mu$ on $X$ in the above sense is equivalent to its quasi-invariance
with respect to the Haar system $\l$, so that $(\G,\l,\mu)$ is a measured
groupoid.

A \emph{random walk along the classes of the equivalence relation} $R$
\cite{Kaimanovich98} is determined by a family of probability measures
$\pi_x,\,x\in X$ concentrated on the corresponding classes $[x]$, which is
Borel in the sense that the function $(x,y)\mapsto\pi_x(y)$ on $R$ is Borel.
The associated invariant Markov operator on the groupoid $\G$ has the
transition probabilities
$$
p\bigl( (x,y), (x,z) \bigr) = \pi_y(z) \;,
$$
which are obviously absolutely continuous with respect to the Haar system.

The same construction is applicable to the leafwise diffusion processes on
\emph{measured Riemannian foliations} as well. In this case the role of the
Haar system can be played by the leafwise Riemannian volumes, see
\cite{Kaimanovich-Lyubich01}.

\medskip

The next 3 examples are obtained from random walks on groups by further
``randomizing'' them in various ways.

\medskip

\num{iii} \emph{Random walks with internal degrees of freedom (RWIDF).} In this
model first introduced by Kr\'{a}mli and Sz\'{a}sz \cite{Kramli-Szasz83} the random
walk on a countable group $G$ is driven by a Markov chain on a space $X$ which
describes the internal or hidden ``degrees of freedom'' of the observed process
on $G$. The state space of the RWIDF is the product $\wt X=X\times G$, and the
transition probabilities $\pi_{(x,g)}=g\pi_x$ are equivariant with respect to
the action of the group on itself, where $\pi_x$ are probability measures on
$\wt X$ indexed by points $x\in X$. In particular, when $X$ is a singleton this
is just the usual random walk on $G$. In the terminology from
\cite{Kaimanovich95} the Markov operator of RWIDF is a \emph{covering operator}
with the \emph{deck group} $G$. The class of covering Markov operators (in
other words, of RWIDF) includes, in particular, the Brownian motion on covering
Riemannian manifolds.

Let us define the corresponding groupoid $\G$ as the product of the groupoid
associated with the full equivalence relation on the set $X$ and the groupoid
associated with the group $G$ (see \secref{sec:examples}). Namely, its set of
morphisms and set of objects are
\begin{equation} \label{eq:big}
\G = X \times X \times G \;, \qquad \Go=X \;,
\end{equation}
respectively, with obvious definitions of the structure maps. Then the random
walk with internal degrees of freedom gives rise to an invariant Markov
operator $P$ on the groupoid $\G$, for which the corresponding probability
measures on the fibers $\G^x$ of the target map (see \propref{prop:extension})
are precisely the measures $\pi_x$.

Any measure $\mu$ on $X=\Go$ is obviously quasi-invariant with respect to the
Haar system $\l$ on $G$ consisting of products of $\mu$ and the counting
measure $\l_G$ on $G$, which gives a structure of a measured groupoid on $\G$.
Absolute continuity of the transition probabilities of the operator $P$ is then
equivalent to absolute continuity of the measures $\pi_x$ with respect to the
product $\mu\otimes \l_G$, or just to absolute continuity of the transition
probabilities of the quotient chain on $X$ with respect to the measure $\mu$.

Under suitable assumptions one can pass in the above construction from
countable groups $G$ to general locally compact groups. However, the situation
when the action of $G$ on the state space $\wt X$ is not free does not readily
fit into this scheme. One has either to pass to the more general construction
of invariant operators on homogeneous spaces of groupoids or to lift the Markov
operator from the homogeneous space to an appropriately defined invariant
operator on the groupoid, see \cite{Kaimanovich-Woess02},
\cite{SaloffCoste-Woess02}

\medskip

\num{iv} \emph{Random walks with random transition probabilities (RWRTP).} This
is a specialization of the above model of RWIDF to the situation when the
quotient chain on $X$ is deterministic, i.e., consists in moving along the
orbits of a certain transformation $T$. Then for any $x\in X$ the corresponding
transition probability $\pi_x$ on $X\times G$ is concentrated on the set
$\{Tx\}\times G$, i.e., can be considered as a probability measure on $G$. Let
us endow the space $X$ with a $T$-invariant probability measure $\mu$. The
associated RWRTP consists in picking up randomly (with distribution $\mu$) a
point $x\in X$ and performing the first step of the random walk on $G$ with the
jump distribution $\pi_x$, then at the next moment of time passing to the point
$Tx$ and performing the next step on $G$ with the jump distribution $\pi_{Tx}$,
etc. Therefore, the observed ``random'' Markov chain on $G$ is homogeneous in
space, but not in time, its transition probabilities at time $n$ being
$p_n(g,gh)=\pi_{T^n x}(h)$, see \cite{Kaimanovich-Kifer-Rubshtein02}.

In order to make the transition probabilities of the invariant Markov operator
arising in this model absolutely continuous one has to pass from the ``big''
groupoid $\G$ \eqref{eq:big} defined in the previous example to the subgroupoid
$\G'\subset\G$ generated by the supports of the transition probabilities. In
this case $\G'$ is the product of the groupoid of the orbit equivalence
relation $R_T$ of the transformation $T$ and the group $G$. Its set of
morphisms is
$$
\G' = \{(x, T^k x, g): x\in X, g\in G, k\in\Z \} \subset \G \;.
$$
and the set objects is the same as for $\G$, i.e., $X$. Then the products of
counting measures on the classes of $R_T$ and the counting measure on the group
$G$ will provide a Haar system on $\G'$, and the transition probabilities of
the arising invariant Markov operator on $\G'$ will obviously be absolutely
continuous with respect to this Haar system.

\medskip

\num{v} \emph{Random walks in random environment (RWRE).} This (historically,
the first model of a randomization of the usual random walk on a group) is yet
another specialization of RWIDF in a sense opposite to RWRTP: the arising
random Markov chains on the group $G$ are homogeneous in time, but not in
space. It was introduced by Solomon \cite{Solomon75} and profoundly studied
later (mostly for abelian groups with few exceptions, though).

Denote by $\Mr(G)$ the space of (infinite) configurations $\om=\{\om_g\}_{g\in
G}$ on $G$ with the values in the space $\Pr(G)$ of probability measures on
$G$. The group $G$ acts on $\Mr(G)$ by translations
$(g\om)_{g'}=\om_{g^{-1}g'}$. The space $\Mr(G)$ can be identified with the
space of all Markov operators on $G$: the transition probabilities of the
Markov operator $P_\om$ determined by $\om$ are
$$
p_\om(g,gh)=\om_g(h)=(g^{-1}\om)_e(h)=p_{g^{-1}\om}(e,h) \;,
$$
i.e., $\om_g$ is the distribution of the ``right increment'' $h$ at the point
$g$.

Let now $X$ be a $G$-space with a $G$-quasi-invariant measure $\mu$ (the
\emph{space of environments}) endowed with a map $x\mapsto\pi_x\in\Pr(G)$. By
equivariance this map can be extended to a map from $X$ to $\Mr(G)$: we shall
consider $\pi_x$ as the value of the configuration $\om(x)$ at the group
identity $e$ and then put
$$
\om(x) = \{\pi_{g^{-1}x}\}_{g\in G} \in \Mr(G) \;.
$$
The associated RWRE consists then in picking up randomly (with distribution
$\mu$) a point $x\in X$ and running on $G$ the ``random'' Markov chain with the
transition probabilities $p_x(g,gh)=\pi_{g^{-1}x}(h)$ (actually, instead of the
space $(X,\mu)$ one could consider directly the space $\Mr(G)$ endowed with the
image of the measure $\mu$ under the map $x\mapsto\om(x)$).

This model gives rise to the covering Markov chain (RWIDF) on the space $\wt
X=X\times G$ with the transition probabilities $\wt\pi_{(x,g)}=g\wt\pi_x$,
where $\wt\pi_x$ are the images of the measures $\pi_x$ under the map $h\mapsto
(h^{-1}x,h)$. Sample paths $(x_n,g_n)$ of this chain have the following
interpretation: $(g_n)$ is a sample path of the Markov chain run in the
environment $g_0^{-1}x_0$, whereas $(x_n)$ is a sample path of the Markov chain
in the space of environments corresponding to the so-called ``moving coordinate
system'' (in which the Markov particle is always situated at the group
identity, whereas the environment around it changes).

In the same way as in the case of RWRTP above, in order to make the transition
probabilities of the arising invariant Markov operator absolutely continuous
one has to pass from the groupoid $\G$ \eqref{eq:big} to a subgroupoid
$\G'\subset\G$ generated by the supports of transition probabilities. In the
case of RWRE
$$
\G'= \{(x,g^{-1}x,g): x\in X, g\in G\} \;,
$$
i.e., this is precisely the groupoid determined by the action of $G$ on $X$.

\section{The Liouville property} \label{sec:liouville}

\subsection{Bounded harmonic functions and the Liouville property}

The classical Liouville theorem asserts absence of bounded harmonic functions
on the Euclidean space. The notion of a harmonic function (based on the mean
value property) can in fact be defined for an arbitrary Markov chain. Namely, a
Borel function $f$ on the state space $X$ of a Markov chain is called harmonic
if $f=Pf$, where $P$ is the transition operator of the Markov chain. In the
measure theoretical setup, given a $P$-adapted measure $m$ on $X$, the classes
(mod 0) of $P$-harmonic functions form a closed subspace
$H^\infty(X,m,P)\subset L^\infty(X,m)$. The operator $P$ is then called
\emph{Liouville} (with respect to the measure $m$) if the space
$H^\infty(X,m,P)$ consists of constant functions only.

The link between the Liouville property and amenability is based on the
so-called \emph{\mbox{0--2} laws} for Markov operators due to Derriennic
\cite{Derriennic76} (also see \cite{Kaimanovich92}), which assert that absence
of non-constant bounded harmonic functions is equivalent to asymptotic
independence of $n$-step transition probabilities of initial states. This is
precisely what is needed for constructing approximatively invariant sequences
of probability measures from condition \eqref{eq:def-amen}. Yet another, less
constructive, way of connecting the Liouville property with amenability
consists in the observation that any fixed \emph{measure-linear mean} on $\Z_+$
when applied to the values of an arbitrary bounded measurable function on the
state space of a Liouville Markov operator along the sample paths of the
associated Markov chain provides a projection onto the space of constants which
is invariant with respect to all the symmetries of the operator
\cite{Connes-Feldman-Weiss81}, \cite{Lyons-Sullivan84},
\cite{Kaimanovich-Fisher98}.

In the case of an invariant Markov operator $P$ on a measured groupoid
$(\G,\l,\mu)$ any measured function on $\G$ which is constant on a.e. fiber
$\G^x$ of the target map $\t:\G\to\Go$ is necessarily harmonic by formula
\eqref{eq:target preserved}, so that it makes sense to talk about the Liouville
property for each of the operators $P_x$ on the fibers $\G^x$.

\subsection{Fiberwise Liouville operators}

\begin{defn} \label{def:fiberwise Liouville operator}
An invariant Markov operator $P$ on a measured groupoid $(\G,\l,\mu)$ is called
\emph{fiberwise Liouville} if for $\mu$-a.e. $x\in\Go$ the operator
$P_x:L^\infty(\G^x,\l^x)\toitself$ is Liouville. A measured groupoid
$(\G,\l,\mu)$ is called \emph{Liouville} if it carries a fiberwise Liouville
invariant Markov operator.
\end{defn}

\begin{thm} \label{thm:liouville measured}
Any Liouville measured groupoid is amenable.
\end{thm}

\begin{proof}
By one of the 0--2 laws a Markov operator $P:L^\infty(X,m)\toitself$ is
Liouville if and only if for any two probability measures $\th_1,\th_2\prec m$
\begin{equation} \label{eq:0-2}
\left\| \frac1{n+1}\sum_{k=0}^n (\th_1-\th_2)P^k \right\| \toto_{n\to\infty} 0
\;.
\end{equation}
Note that if the transition probabilities $\pi_x$ of the operator $P$ are
absolutely continuous with respect to the measure $m$, then it is sufficient to
consider in formula \eqref{eq:0-2} just the $\d$-measures
$\th_i=\d_{x_i},\,x_i\in X\;(i=1,2)$.

Let now $P:L^\infty(\G,\l\star\mu)\toitself$ be a fiberwise Liouville invariant
Markov operator. Take a measurable system of absolutely continuous probability
measures $\th^x\prec\l^x$ on the fibers $\G^x$ of the target map, and let
$$
\th_n^x = \frac1{n+1}\sum_{k=0}^n \th^x P^k \prec \l^x \;.
$$
By $\G$-invariance of the operator $P$ for any $\g\in\G$
\begin{equation} \label{eq:difference}
\left\| \g\th_n^{\s(\g)} - \th_n^{\t(\g)} \right\| =
 \left\| \frac1{n+1}\sum_{k=0}^n \left( \g\th^{\s(\g)} - \th^{\t(\g)} \right) P^k
 \right\|\;,
\end{equation}
where $\g\th^{\s(\g)}$ and $\th^{\t(\g)}$ are probability measures on the fiber
$\G^{\t(\g)}$ absolutely continuous with respect to $\l^{\t(\g)}$. Since the
fiberwise operators are a.e. Liouville, the 0--2 law \eqref{eq:0-2} implies
a.e. convergence of \eqref{eq:difference} to zero, and therefore, in view of
formula \eqref{eq:def-amen}, amenability of the groupoid $(\G,\l,\mu)$.
\end{proof}

\begin{rem}
Another (non-constructive) proof of \thmref{thm:liouville measured} can be
given by using a measure-linear mean along the sample paths of the Markov chain
to obtain an invariant mean from $L^\infty(\G)$ to $L^\infty(\Go)$, see the
proof of \thmref{thm:extension} below.
\end{rem}

\subsection{Applications and examples}

Particular cases of \thmref{thm:liouville measured} were earlier established
for groups \cite{Azencott70}, \cite{Furstenberg73}, equivalence relations and
foliations \cite{Connes-Feldman-Weiss81}, for the Brownian motion on covering
manifolds \cite{Lyons-Sullivan84}, for general covering Markov operators
\cite{Kaimanovich95} as well as for various models of randomization of the
usual random walk on a discrete group, see \cite{Sunyach87},
\cite{Kaimanovich-Kifer-Rubshtein02}.

It is plausible that the converse may also be true:

\begin{conj} \label{conj:liouville}
Any amenable measured groupoid is Liouville.
\end{conj}

This is known to be the case for groups \cite{Rosenblatt81},
\cite{Kaimanovich-Vershik83} (any amenable group carries a random walk with the
trivial Poisson boundary; it had been previously conjectured by Furstenberg
\cite{Furstenberg73}). The proof of \conjref{conj:liouville} in full generality
should presumably follow the same strategy of constructing a Liouville operator
from F{\o}lner sets on the groupoid. Other known particular cases are the
groupoids associated with discrete equivalence relations (in view of the
Connes--Feldman--Weiss theorem it is the orbit equivalence relation of a
$\Z$-action \cite{Connes-Feldman-Weiss81}) and with group actions (in a
somewhat weaker form, though; by \cite{Elliott-Giordano93},
\cite{Adams-Elliott-Giordano94} any amenable measure class preserving action of
a locally compact group $G$ can be realized as the action on the Poisson
boundary of an appropriate $G$-invariant operator, see below
\secref{sec:examples poisson}; proving \conjref{conj:liouville} for group
actions would provide an alternative link between amenability and the Poisson
boundary). Note that the proof of hyperfiniteness of amenable equivalence
relations in \cite{Connes-Feldman-Weiss81} (also see \cite{Kaimanovich97}) is
non-constructive and uses the Zorn Lemma; looking for a more direct F{\o}lner type
argument would provide an insight into the general case.

\begin{rem}
There is a generalization to group extensions of the aforementioned existence
of a Liouville random walk on any amenable locally group. Namely, if $G$ is a
locally compact group, and $H$ its closed normal subgroup, then $H$ is amenable
if and only if for any Borel probability measure $\pi'$ on the quotient group
$G'=G/H$ there exists a Borel lift $\pi$ to $G$ such that the Poisson
boundaries of the random walks $(G,\pi)$ and $(G',\pi')$ are canonically
isomorphic \cite{Kaimanovich02a}. In spite of having the same spirit as
\conjref{conj:liouville}, this result does not seem to have an obvious
interpretation in groupoid terms.
\end{rem}

\subsection{Amenability of group actions}

We shall now give a specialization of \thmref{thm:liouville measured} to the
particular case of (groupoids associated with) group actions. In this situation
the amenability of a measure class preserving action of a locally compact group
$G$ on a measure space $(X,\mu)$ is equivalent to existence of a sequence of
measurable maps $\th_n$ from the action space $X$ to the space $\Pr(G)$ of
probability measures on $G$ which are \emph{approximatively equivariant} in the
sense that $\|g\th_n(x)-\th_n(gx)\|\to 0$ weakly, cf. formula
\eqref{eq:def-amen}. Recall that a Borel $G$-space $S$ is called proper if it
carries a Borel $G$-invariant system of probability measures on $G$
\cite[Definition 2.1.2]{Anantharaman-Renault00} (for continuous actions on
locally compact spaces it follows prom properness in the usual sense). Such a
system allows one to lift any probability measure from $S$ to $G$ in an
equivariant measurable way. We shall need this property in the measure
theoretical setup and say that an action of a locally compact group $G$ on a
measure space $(S,m)$ is \emph{proper} if there exists a measurable equivariant
map from the space of probability measures on $S$ absolutely continuous with
respect to $m$ to the space of probability measures on $G$.

\begin{thm} \label{thm:amenab actions-measurable}
Let $G$ be a locally compact group acting measure preserving and properly on a
measure space $(S,m)$, and let $P:L^\infty(S,m)\toitself$ be a $G$-invariant
Markov operator. Suppose that the group $G$ also has a measure class preserving
action on another space $(X,\mu)$, and there is a measurable $G$-equivariant
map assigning to points $x\in X$ projective classes of minimal positive Borel
$P$-harmonic functions $\f_x$. Then the action of $G$ on $(X,\mu)$ is amenable.
\end{thm}

\begin{proof}
Recall that any positive $P$-harmonic function $\f$ determines a new Markov
operator $P^\f f=P(\f f)/\f$ on $L^\infty(S,m)$ called the \emph{Doob
transform} of the operator $P$ determined by the function $\f$ \cite{Revuz84}.
In other words, the transition probabilities of the Doob transform are
determined by the formula
\begin{equation} \label{eq:doob}
\frac{d\pi^\f_s}{d\pi_s}(t)= \frac{\f(t)}{\f(s)} \;,
\end{equation}
where $\pi_s$ are the transition probabilities of the operator $P$. The Doob
transform remains the same if the function $\f$ is multiplied by a positive
constant, so that actually it depends just on the projective class of $\f$, or,
in other words, on the multiplicative cocycle $(s,t)\mapsto\f(t)/\f(s)$. One
can easily see that the minimality of the function $\f$ is equivalent to the
Liouville property for the operator $P^\f$, which provides the required link
between minimality and the Liouville property.

Strictly speaking, the situation considered in \thmref{thm:amenab
actions-measurable} is slightly different from the setup of
\thmref{thm:liouville measured} as here we arrive at an invariant Markov
operator on a homogeneous space of the action groupoid rather than just on the
groupoid itself (cf. the discussion at the end of Example (iii) in
\secref{sec:examples of operators}), which is why we shall briefly outline the
rest of the proof.

Let us fix a probability measure $\th\prec m$ on $S$, and consider the family
$$
\th_n^x = \frac1{n+1} \sum_{k=0}^n \th (P^x)^k \prec m
$$
of probability measures on $S$ parameterized by points $x\in X$, where
$P^x=P^{\f_x}$ are the Doob transforms associated with the functions $\f_x$.
Then, in the same way as in the proof of \thmref{thm:liouville measured}, the
systems $\th_n^x$ are approximatively equivariant, i.e., for any $g\in G$
$$
\| g\th_n^x - \th_n^{gx} \| \to 0 \;.
$$
Since the action of $G$ on $(S,m)$ is proper, the measures $\th_n^x$ can now be
equivariantly and measurably lifted from $S$ to $G$ to provide an
approximatively equivariant sequence of measurable maps from $X$ to $\Pr(G)$.
\end{proof}

\begin{rem}
\thmref{thm:amenab actions-measurable} and its topological analogue
\thmref{thm:amenab actions} carry over \emph{verbatim} to the situation when
$\f_x$ are \emph{$\l$-harmonic} minimal functions for a certain fixed
eigenvalue $\l>0$ (i.e., $P\f_x=\l\f_x$). In this case the definition of the
Doob transform has to be modified by dividing the right-hand side of formula
\eqref{eq:doob} by $\l$.
\end{rem}

\begin{rem}
The metric characterization of minimal harmonic functions resulting from
applying the 0--2 law to the corresponding Doob transform first appeared in
author's paper \cite{Kaimanovich83}.
\end{rem}

\begin{rem}
The properness assumption in \thmref{thm:amenab actions-measurable} is
essential. In a sense, it says that the measurable structures on the group $G$
and on the $G$-space $S$ agree. For instance, take a free dense subgroup $F$ of
a compact group $K$, and consider on $K$ the random walk determined by a
transition probability $\pi$ supported by the generating set of $F$. Then the
Poisson boundary of the associated transition operator on the space $(K,\l)$
(where $\l$ is the Haar measure on $K$) is trivial (there are no measurable
$\pi$-harmonic functions on $K$, which follows from the ergodicity of the
action of $F$ and the fact that $\l$ is a finite stationary measure of this
random walk), but the action of the free group on a singleton is not amenable.
This example illustrates the importance of the choice of an ambient measurable
structure in the definition of the Poisson boundary.
\end{rem}

\begin{rem}
In the Borel setup, when the map $x\mapsto\f_x$ is well-defined for \emph{all}
points $x\in X$, the argument from the proof of \thmref{thm:amenab
actions-measurable} is applicable to an arbitrary quasi-invariant measure $\mu$
on $X$ to provide the \emph{measurewise amenability}
\cite{Anantharaman-Renault00} (other terms: \emph{universal amenability}
\cite{Adams96}, \emph{measure-amenability} \cite{Jackson-Kechris-Louveau02}) of
the action of $G$ on $X$; cf. \remref{rem:measurewise} below).
\end{rem}

See \secref{sec:applications topological} for examples of application of
\thmref{thm:amenab actions-measurable} (or, rather, of its topological
refinement \thmref{thm:amenab actions}).

\section{The Poisson extension} \label{sec:poisson extension}

\subsection{The Poisson boundary} \label{sec:poisson boundary}

Let $m$ be a $P$-adapted Borel $\si$-finite measure of a Markov operator $P$ on
a state space $X$, so that the operator $P$ acts on the space $L^\infty(X,m)$.
We shall assume that $(X,m)$ is a \emph{Lebesgue space} (in particular, this is
always the case when $X$ is \emph{Polish}). Then the associated path space
$(X^{\Z_+},\P_m)$ is also a Lebesgue space, and the space $\Ga=\Ga(X,m,P)$ of
the ergodic components of the time shift $T$ in the path space is called the
\emph{Poisson boundary} of the operator $P$. By definition, there is a
canonical measurable projection $\bnd$ from the path space onto the Poisson
boundary constant along the orbits of the time shift, and the Poisson boundary
is the maximal quotient of the path space with this property. The Poisson
boundary (which is defined in the measure theoretical category only!) is
endowed with the \emph{harmonic measure class} $[\nu_m]=\bnd[\P_m]$. For
convenience we shall fix a probability measure $\nu\in [\nu_m]$. For instance,
one can take $\nu=\bnd\P_\th$ for any probability measure $\th$ on the state
space equivalent to $m$. For any initial distribution $\th\prec m$ the
associated \emph{harmonic measure} $\nu_\th=\bnd\P_\th$ is absolutely
continuous with respect to the harmonic measure class. If the operator $P$ has
absolutely continuous transition probabilities then the individual harmonic
measures $\nu_x$ (corresponding to the initial distributions $\d_x,\,x\in X$)
are also well-defined and absolutely continuous with respect to the harmonic
measure class.

The space of (classes mod 0 of) bounded harmonic functions $H^\infty(X,m,P)$ of
the operator $P$ is canonically isomorphic to the space
$L^\infty(\Ga,[\nu_m])$. For Markov operators with absolutely continuous
transition probabilities this isomorphism is established by the \emph{Poisson
formula}
$$
f(x) = \langle \nu_x, \wh f \rangle = \int \wh f(\wh\ga) \,
\Pi(x,\ga)\,d\nu(\ga) \;,
$$
where $\Pi(x,\ga)=d\nu_x/d\nu(\ga)$ is the \emph{Poisson kernel}. For
$[\nu_m]$-a.e. point $\ga\in\Ga$ the function $\Pi(\cdot,\ga)$ is a minimal
$P$-harmonic function. Actually, the Poisson formula can be given sense in full
generality as well by defining the individual harmonic measures (not
necessarily absolutely continuous with respect to the harmonic measure class
anymore!) by using Rokhlin's theorem on conditional probabilities in Lebesgue
spaces (cf. the proof of \thmref{thm:extension} below).

See \cite{Kaimanovich92} and the references therein for more detailed
information on the Poisson boundary of Markov operators.

\subsection{The Poisson boundary of invariant operators}

Let now $P$ be an invariant Markov operator on a measured groupoid
$(\G,\l,\mu)$, so that the measure $\l\star\mu$ is $P$-adapted. We denote
sample paths from $\G^{\Z_+}$ by $\ov\g=(\g_0,\g_1,\dots)$. Since the target
map $\t$ is constant along the sample paths of an invariant Markov operator by
formula \eqref{eq:target preserved}, it can be extended to a projection map
(also denoted $\t$) from the path space to $\Go$. Therefore, the action of $\G$
on itself extends to a coordinate-wise action of $\G$ on the path space by the
formula
$$
\g\ov\g=(\g\g_0,\g\g_1,\dots) \;, \qquad \s(\g)=\t(\ov\g) \;,
$$
where $\g_n$ are the coordinates of the sample path $\ov\g$. Since the measure
$\l\star\mu$ on $\G$ is quasi-invariant with respect to $(\G,\l)$ (see
\secref{sec:measured groupoids}), the associated measure $\P_{\l\star\mu}$ on
the path space is also $(\G,\l)$-quasi-invariant. Further, since the action of
the time shift on the path space commutes with the action of $\G$, we obtain
that this action descends to the Poisson boundary $\Ga$ of the operator $P$,
and that the harmonic measure class on $\Ga$ is quasi-invariant with respect to
this action. The target map descends to $\Ga$ from the path space; its fibers
are the Poisson boundaries of the fiberwise Markov operators
$P_x:L^\infty(\G^x,\l^x)\toitself$ (cf.
\cite[Proposition~1.11]{Kaimanovich-Kifer-Rubshtein02}).

\vfill

\pagebreak

\subsection{The Poisson extension}

\begin{defn}
Given an invariant Markov operator $P$ on a measured groupoid $(\G,\l,\mu)$, we
shall call the \emph{Poisson extension} the measured groupoid $\wt\G=\G\lsd\Ga$
associated with the action of $\G$ on the Poisson boundary $\Ga$ of the
operator $P$.
\end{defn}

The following result is a generalization of \thmref{thm:liouville measured} (if
an invariant operator is fiberwise Liouville, then its Poisson extension is
just the original groupoid).

\begin{thm} \label{thm:extension}
The Poisson extension of any invariant Markov operator on a measured groupoid
is amenable.
\end{thm}

\begin{proof}
One way of proving \thmref{thm:extension} is to deduce it from
\thmref{thm:liouville measured} by showing that the groupoid $\wt\G$ carries a
natural fiberwise Liouville invariant Markov operator. Such an operator is
obtained by conditioning the original operator $P$ by the points of the Poisson
boundary. This is easily done in the case when $P$ has absolutely continuous
transition probabilities. Namely, a.e. point $\ga\in\Ga$ determines the
conditional Markov operator $P^\ga$ on the fiber $\G^{\t(\ga)}$, for which the
measure $\l^{\t(\ga)}$ is adapted. These are the Doob transforms corresponding
to the fiberwise Poisson kernels. In other words, the transition probabilities
of $P^\ga$ satisfy the relation \eqref{eq:doob}
$$
\frac{d\pi_\g^\ga}{d\pi_\g}(\g') = \frac{d\nu_{\g'}}{d\nu_\g}(\ga) \;,
$$
where $\nu_\g$ are the harmonic measures on the Poisson boundary (it is here
that we need the absolute continuity of the transition probabilities which
guarantees existence of fiberwise Poisson kernels). Since the Poisson kernel
consists of minimal harmonic functions, the operators
$P^\ga:L^\infty\left(\G^{\t(\ga)},\l^{\t(\ga)}\right)\toitself$ are Liouville.

It remains to notice that the conditional operators $P^\ga$ can be interpreted
as fiberwise operators of an invariant Markov operator $\wt P$ on the Poisson
extension $\wt\G$ and to apply \thmref{thm:liouville measured}. Indeed, the
elements of $\wt\G=\G\lsd\Ga$ are the triples $(\ga,\g,\g^{-1}\ga)$ (see
\secref{sec:homogeneous spaces}). Fixing the target $\ga=\t(\ga,\g,\g^{-1}\ga)$
(which corresponds to conditioning by $\ga$ as an element of the Poisson
boundary) we may identify the corresponding fiber $\wt\G^\ga$ with
$\G^{\t(\ga)}$ by the formula $(\ga,\g,\g^{-1}\ga)\leftrightarrow\g$ and then
define the transition probabilities
$$
\wt\pi_{(\ga,\g,\g^{-1}\ga)} = \pi_\g^\ga \;,
$$
which are clearly $\wt\G$-invariant in view of $\G$-invariance of the
probabilities $\pi_\g$, cf. Example~(v) from \secref{sec:examples of
operators}.

\medskip

Let us sketch how this argument can be carried over to the general case. Since
the measure $\l\star\mu$ is $P$-adapted, for the associated measure in the path
space $T\P_{\l\star\mu}\prec\P_{\l\star\mu}$ (here $T$ is the time shift). If
$\th$ is a probability measure equivalent to $\l\star\mu$, then clearly
\begin{equation} \label{eq:prec}
T\P_\th\prec\P_\th \;.
\end{equation}
Denote by $\P_\th^\ga,\,\ga\in\Ga$ the $T$-ergodic components of the measure
$\P_\th$, i.e., its conditional measures with respect to the Poisson boundary.
By definition of the Poisson boundary these measures are Markov and, by
\eqref{eq:prec}, $T\P_\th^\ga\prec\P_\th^\ga$ for a.e. $\ga\in\Ga$. Therefore
the one-dimensional distributions $\th^\ga$ of the measures $\P_\th^\ga$ at
time 0 are adapted with respect to the corresponding conditional Markov
operators $P^\ga$ (although these measures may well be singular with respect to
the Haar measures $\l^{\t(\ga)}$). The rest of the argument then goes in the
same way as in the absolutely continuous case.

\medskip

However, it is more convenient to use for proving \thmref{thm:extension}
another approach based on using \emph{measure-linear means}. Mokobodzki (see
\cite{Meyer73}, \cite{Fisher87}) proved that there exists an invariant mean
$\xi$ on $\Z_+$ (called \emph{measure-linear} or \emph{medial}) with the
following remarkable property: it is \emph{universally measurable} as a map
from the product space $[-1,1]^{\Z_+}$ to $[-1,1]$, i.e., the integral in the
right hand side below is well-defined for any Borel probability measure $\mu$
on $[-1,1]^{\Z_+}$, and
$$
\xi\left\{\int\ab\,d\mu(\ab)\right\}=\int\xi\{\ab\}\,d\mu(\ab) \;.
$$

In view of the definition from \secref{sec:amenable groupoids} for proving
amenability of the groupoid $\wt\G$ we have to construct an invariant mean
$\PPi:L^\infty(\wt\G)\to L^\infty(\Ga)$ (recall that the space of objects of
$\wt\G=\G\lsd\Ga$ is the Poisson boundary $\Ga$). As above we shall
parameterize $\wt\G$ by the map
\begin{equation} \label{eq:coordinates}
(\g,\ga)\leftrightarrow (\ga,\g,\g^{-1}\ga) \;, \qquad \g\in\G,\,\ga\in\Ga \;,
\end{equation}
where $\t(\g)=\t(\ga)$. Then the target map $\wt\G\to\wt\G^{0}$ is just
$(\g,\ga)\mapsto\ga$, and the measure class on $\wt\G$ with respect to which we
consider the space $L^\infty(\wt\G)$ is $d\l^{\t(\ga)}(\g)\,d[\nu](\ga)$, where
$[\nu]$ is the harmonic measure class on $\Ga$ corresponding to the initial
distribution $\l\star\mu$ on $\G$. In the coordinates \eqref{eq:coordinates}
the left action of $\wt\G$ on itself coincides just with the diagonal action of
$\G$
$$
\h(\g,\ga) \mapsto (\h\g,\h\ga) \;, \qquad \s(\h)=\t(\g)=\t(\ga) \;.
$$

Let us fix a reference system of probability measures $\rho=\{\rho^x\}$ on the
fibers of the target map of $\G$ equivalent to the Haar system $\l$ and put for
any $F\in L^\infty(\wt\G)$ and any sample path $\ov\g=(\g_n)$ on $\G$
$$
\ov\PPi F(\ov\g) = \xi \left\{ \int F(\g_n\h,\bnd\ov\g) \,d\rho^{\s(\g_n)} (\h)
\right\} \;.
$$
Note that \emph{a priori} the images of the measure $\P_{\l\star\mu}$ under the
maps $\ov\g\mapsto (\g_n,\bnd\ov\g)$ may well be singular with respect to the
quasi-invariant measure class on $\wt\G$ (we keep using the coordinates
\eqref{eq:coordinates}), however the additional integration with respect to the
measures $\rho^x$ guarantees that $\ov\PPi F$ is well-defined as an element of
$L^\infty(\G^{\Z_+},\P_{\l\star\mu})$. Since $\xi$ is a mean, the function
$\ov\PPi F$ is shift invariant, so that it descends to a measurable function
$$
\PPi F(\bnd\ov\g) = \ov\PPi F(\ov\g)
$$
on the Poisson boundary. Since $\xi$ is measure-linear, $\PPi$ is a mean from
$L^\infty(\wt\G)$ to $L^\infty(\Ga)\cong L^\infty(\wt\G^{0})$. Finally,
$\wt\G$-invariance of $\PPi$ obviously follows from its definition.
\end{proof}

\subsection{Examples} \label{sec:examples poisson}

\thmref{thm:extension} has been earlier proved in the following particular
cases.

For the groupoid $\G$ associated with a locally compact group $G$ and the
invariant Markov operator determined by a random walk $(G,\pi)$ (see
\secref{sec:examples of operators}) amenability of the Poisson extension
amounts to ergodicity of the action of $G$ on the Poisson boundary of the
associated Markov operator on the space $L^\infty(G,\l)$ (where $\l$ is the
Haar measure on $G$). It was proved by Zimmer \cite{Zimmer78} for an arbitrary
measure $\pi$ on $G$. Note that the formulation of Theorem 5.1 in
\cite{Zimmer78} contains the requirement that the measure $\pi$ be
\emph{\'{e}tal\'{e}e} or \emph{spread out} (i.e., have a convolution power non-singular
with respect to the Haar measure). The reason is that there he used a
definition of the Poisson boundary which differs from ours; it was supposed to
represent \emph{all} bounded harmonic functions rather then just classes mod 0
from $L^\infty(G,\l)$. However, Zimmer essentially deals precisely with the
Poisson boundary in our sense, and proves its amenability without any further
assumptions on the measure $\pi$ in his Theorem 5.2.

Zimmer used the fact that the increments of a random walk on a group form a
stationary sequence. As it was pointed out in \cite{Connes-Woods89}, his method
does not seem to work in the non-stationary case. Connes and Woods
\cite{Connes-Woods89} proved amenability of the action of a locally compact
group $G$ on the Poisson boundary for the so-called \emph{matrix-valued random
walks} on $G$ which are $G$-invariant Markov chains on the product of $G$ by a
certain countable set subject to some additional assumptions on transition
probabilities. Note that actually what they call the Poisson boundary is rather
the \emph{tail boundary} (the quotient of the path space by the
\emph{synchronous} asymptotic equivalence relation). Although it coincides with
our Poisson boundary (the quotient of the path space by the \emph{asynchronous}
asymptotic equivalence relation) for matrix-valued random walks in the sense of
\cite{Connes-Woods89}, in general they may differ (the Poisson boundary being a
quotient of the tail boundary, so that amenability of the action on the former
is stronger than on the latter), for instance, see \cite{Kaimanovich92},
\cite{Jaworski95}.

Elliott and Giordano \cite{Elliott-Giordano93} for discrete groups and Adams,
Elliott, Giordano \cite{Adams-Elliott-Giordano94} in the general case later
proved that in fact any measure class preserving amenable action of a second
countable locally compact group can be presented as its action on the Poisson
boundary of an appropriately defined matrix-valued random walk.

\section{Topological amenability} \label{sec:topological}

\subsection{Topological Liouville property}

\thmref{thm:liouville measured} has an analogue in the topological category
(see \secref{sec:top amenab} for a definition of topological amenability) which
can be proved along the same lines by using the 0--2 law. When working in the
topological setup we have to modify the definition of Liouville operators by
saying that an invariant Markov operator is \emph{topologically fiberwise
Liouville} if \emph{all} (rather than almost all as in
\defref{def:fiberwise Liouville operator}) fiberwise operators
$P_x:L^\infty(\G^x,\l^x)\toitself$ are Liouville.

\begin{thm} \label{thm:top amenab}
Let $\G$ be a locally compact topological groupoid with a continuous Haar
system, and let $P$ be an invariant Markov operator on $\G$ with continuous
densities. If the operator $P$ is topologically fiberwise Liouville, then the
groupoid $\G$ is topologically amenable.
\end{thm}

\begin{proof}
We shall use the fact that the Cesaro averages in the formulation of the 0--2
law for the triviality of the Poisson boundary can be replaced with any
sequence of probability measures on $\Z_+$ strongly convergent to invariant
mean on $\Z$ \cite{Kaimanovich92}. By taking for such a sequence the binomial
distributions on $\Z_+$ the convergence in formula \eqref{eq:0-2} can be made
monotone and therefore uniform on compacts.

More precisely, let us consider the Markov operator $Q=(P+P^2)/2$. The
transition probabilities of the operator $Q$ are the averages of the time 1 and
time 2 transition probabilities of the operator $P$. Obviously, the operator
$Q$ is also invariant, and it is fiberwise Liouville simultaneously with the
operator $P$. For the operator $Q$ the Poisson boundary coincides with the tail
boundary, and therefore the corresponding 0--2 law (see \cite{Kaimanovich92})
implies that
\begin{equation} \label{eq:pp2}
\left\| \left( \d_\g-\d_{\t(\g)} \right) Q^n \right\| \toto_{n\to\infty} 0
\qquad \forall\,\g\in\G \;,
\end{equation}
or, in other words, that
$$
\left\| 2^{-n} \sum_{k=0}^n \binom{n}{k} \left( \d_\g-\d_{\t(\g)} \right)
P^{n+k} \right\| \toto_{n\to\infty} 0 \qquad \forall\,\g\in\G \;.
$$
It is clear that the convergence in formula \eqref{eq:pp2} is monotone. On the
other hand, continuity of densities of the transition probabilities of the
operator $P$ (therefore, of its powers as well) and continuity of the Haar
system implies that the left-hand side of \eqref{eq:pp2} depends on $\g$
continuously. It remains to refer to Dini's theorem, according to which
monotone convergence of continuous functions to a continuous limit on a compact
set is uniform, and to conclude in the same way as in the proof of
\thmref{thm:liouville measured}.
\end{proof}

\subsection{Amenability of group actions}

In the same way as in \secref{sec:liouville}, we shall now give a
specialization of \thmref{thm:top amenab} to the particular case of (groupoids
associated with) group actions. In this situation topological amenability of a
continuous action of a locally compact group $G$ on a locally compact space $X$
is equivalent to existence of a sequence of weak$^*$ continuous maps $\th_n$
from the action space $X$ to the space $\Pr(G)$ of probability measures on $G$
which are \emph{topologically approximatively equivariant} in the sense that
$\|g\th_n(x)-\th_n(gx)\|\to 0$ uniformly on compact subsets of $G\times X$ (cf.
\secref{sec:top amenab}).

\begin{thm} \label{thm:amenab actions}
Let $G$ be a locally compact group acting continuously and properly on a
locally compact space $S$, let $m$ be a $G$-invariant Borel measure on $S$, and
let $P$ be a $G$-invariant Markov operator on $S$ with absolutely continuous
transition probabilities and continuous densities with respect to the measure
$m$. Suppose that the group $G$ also acts continuously on another locally
compact space $X$, and there is a continuous $G$-equivariant map assigning to
points $x\in X$ projective classes of minimal positive $P$-harmonic functions
$\f_x$. Then the action of $G$ on $X$ is topologically amenable.
\end{thm}

\begin{proof}
As in the proof of \thmref{thm:top amenab}, let us pass from the operators
$P_x=P_{\f_x}$ to the operators $Q_x=(P_x+P_x^2)/2$, and consider on $S$ the
probability measures
$$
\th_n(s,x) = \d_s  Q_x^n \;, \qquad s\in X,\; x\in X \;.
$$
The maps $(s,x)\to\th_n(s,x)$ are $G$-equivariant, and
$$
\| \th_n(s,x) - \th_n(s',x) \| \toto_{n\to\infty} 0
$$
uniformly on compact sets. Then for any fixed point $o\in S$
$$
\| g\th_n(o,x) - \th_n(o,gx) \| = \| \th_n(go,gx) - \th_n(o,gx) \|
\toto_{n\to\infty} 0
$$
for any $g\in G$ and $x\in X$ uniformly on compact sets. Since the action of
$G$ on $S$ is proper, the measures $\th_n(o,x)$ can be equivariantly and
continuously lifted from $S$ to $G$ \cite[Corollary 2.1.17]
{Anantharaman-Renault00}, to provide a topologically approximatively
equivariant sequence of maps from $X$ to $\Pr(G)$.
\end{proof}

\begin{rem} \label{rem:measurewise}
In \cite{Biane-Germain02} the following particular cases of \thmref{thm:amenab
actions} were established:
\begin{itemize}
  \item [(i)]
$G$ is a finitely generated group, $S=G$, and $P$ is the Markov operator of a
finitely supported random walk on $G$;
  \item [(ii)]
$G$ is a lattice in an ambient locally compact group $S$, and $P$ is the Markov
operator of an absolutely continuous random walk on $S$.
\end{itemize}
Although the proof in \cite{Biane-Germain02} was also based on using the 0-2
law, it was done in the measure theoretical setting only by applying then the
theorem of Anantharaman-Delaroche and Renault on the equivalence of the
topological amenability and the \emph{measurewise amenability} (i.e., the
amenability of the measured groupoid $(\G,\l,\mu)$ for any measure $\mu$ on
$\Go$ quasi-invariant with respect to a fixed Haar system $\l$) for locally
compact groupoids with a continuous Haar system and countable orbits
\cite[Theorem 3.3.7]{Anantharaman-Renault00}.
\end{rem}

\subsection{Spaces of minimal harmonic functions} \label{sec:applications topological}

As we have already mentioned in \secref{sec:poisson boundary}, the
Radon--Nikodym derivatives of harmonic measures on the Poisson boundary of a
Markov operator with absolutely continuous transition probabilities are minimal
harmonic functions. The \emph{Martin boundary} is a topological counterpart of
the Poisson boundary. Unlike the Poisson boundary (defined as a measure space),
the Martin boundary is a \emph{bona fide} topological space obtained by taking
the closure of the topological state space embedded into the space of functions
on itself via the Green kernel under suitable regularity conditions on
transition probabilities \cite{Revuz84}. The Martin boundary contains all
minimal harmonic functions (the closure of the corresponding subset is
sometimes called the \emph{minimal Martin boundary}). Since the space of
positive harmonic functions is a lattice, any harmonic function can be uniquely
decomposed as an integral of minimal ones. Note that the Martin boundary
considered as a measure space endowed with the representing measures of the
constant function $\1$ coincides with the Poisson boundary. The action of any
symmetry group of the Markov operator extends to the Martin boundary. See
\cite{Kaimanovich96} and the references therein for a discussion of the Martin
boundary for Markov operators on homogeneous spaces.

\thmref{thm:amenab actions} immediately implies:

\begin{thm} \label{thm:martin}
Under conditions of \thmref{thm:amenab actions}, if the subset $M$ of minimal
harmonic functions in the Martin boundary of the operator $P$ is closed, then
the action of $G$ on $M$ is topologically amenable.
\end{thm}

Identification of the space of minimal harmonic functions of a $G$-invariant
Markov operator is, in general, a difficult problem (e.g., see
\cite{Kaimanovich96}). Note, however, that there is no need for the space $X$
from \thmref{thm:amenab actions} to represent \emph{all} minimal harmonic
functions. Appropriate geometrical boundaries were shown to produce in a
continuous way minimal harmonic functions (not necessarily all of them!) in
several situations of hyperbolic flavour:
\begin{itemize}
  \item [(i)]
If $X$ is a \emph{Gromov hyperbolic Riemannian manifold} of bounded geometry
with a spectral gap, $P$ is the Markov operator corresponding to the
\emph{Brownian motion} on $X$, and $\part X$ is the \emph{hyperbolic boundary}
of $X$ \cite{Ancona90}, in particular, if $X$ is a \emph{simply connected
Riemannian manifold with pinched sectional curvatures} and $\part X$ is its
\emph{visibility boundary} \cite{Anderson-Schoen85};
  \item [(ii)]
If $X$ is a \emph{Gromov hyperbolic graph} satisfying the strong isoperimetric
inequality, $P$ is the Markov operator of the \emph{simple random walk} on $X$,
and $\part X$ is the \emph{hyperbolic boundary} of $X$ \cite{Ancona90};
  \item [(iii)]
If $X$ is a \emph{non-compact Riemannian symmetric space}, $P$ is the Markov
operator corresponding to the \emph{Brownian motion} on $X$, and $\part X$ is
the \emph{Furstenberg boundary} of $X$ (it is defined as the space of
\emph{asymptotic classes of Weyl chambers} in $X$, or, equivalently, as the
quotient of the corresponding semi-simple Lie group $G$ by a minimal parabolic
subgroup; for $G=SL(n,\R)$ this is the \emph{flag space} in $\R^n$)
\cite{Furstenberg63}, \cite{Karpelevich65}, see the book
\cite{Guivarch-Ji-Taylor98} for the latest developments, in particular, for a
description of minimal harmonic functions for random walks on general
unimodular groups with Gelfand pairs (Theorem 13.12);
  \item [(iv)]
If $X$ is a locally finite \emph{affine building}, $P$ is the Markov operator
of the \emph{simple random walk} on its set of vertices, and $\part X$ is the
\emph{spherical building at infinity} of $X$ (the space of \emph{asymptotic
classes of sectors}); this case requires more explanations, so that its
discussion is relegated to \secref{sec:buildings} below.
\end{itemize}

In all these cases the operator $P$ agrees with the underlying geometrical
structure on $X$, so that it commutes with the (proper) group of isomorphisms
of $X$. \thmref{thm:amenab actions} then provides topological amenability of
the corresponding boundary actions ``for free'', which gives a unified
generalization of numerous earlier results on amenability of boundary actions
\cite{Bowen77}, \cite{Vershik78}, \cite{Zimmer84}, \cite{Spatzier87},
\cite{Spatzier-Zimmer91}, \cite{Adams94}, \cite{Adams96},
\cite{Robertson-Steger96}, \cite{Ramagge-Robertson96},
\cite{Cutting-Robertson03}:

\begin{thm} \label{thm:geometric amenability}
The action of a closed group of isomorphisms of the space $X$ on the boundary
$\part X$ is topologically amenable in the above cases (i) -- (iv).
\end{thm}

\subsection{Harmonic functions on buildings} \label{sec:buildings}

The set $V$ of vertices of a locally finite affine building $X$ is split into
several types, so that by taking this additional structure into account one can
naturally define several commuting Markov operators (``Laplacians'')
$P_1,P_2,\dots,P_d$ on $V$ (where $d$ is the dimension of the building; see the
references below for details). By the \emph{simple random walk} on $V$ we shall
mean the Markov chain associated with the operator $P=(P_1+P_2+\dots+P_d)/d$.

A function $f$ on $V$ is called \emph{strongly harmonic} if it is harmonic in
the usual sense ($P_i f=f$) for all operators $P_i$. It is known that there is
a continuous map assigning to points $\ga\in\part X$ functions $\f_\ga$ on $V$
which are strongly harmonic and minimal in the cone of non-negative strongly
harmonic functions (minimality follows from the uniqueness of the boundary
decomposition for strongly harmonic functions). More precisely, it follows from
the results of Kato \cite{Kato81} (also see \cite[Theorem
13.12]{Guivarch-Ji-Taylor98}) that this is true for affine buildings associated
with reductive groups over $p$-adic fields, in particular, for all buildings of
dimension at least 3 \cite{Tits86}. Dimension 1 affine buildings are just
trees, whereas for the remaining ``non-classical'' affine buildings of
dimension 2 it was proved by Mantero and Zappa \cite{Mantero-Zappa98},
\cite{Mantero-Zappa00}, \cite{Mantero-Zappa02} (also see \cite{Cartwright99}
for a unified treatment of type $\tilde A_n$ buildings).

Obviously, any strongly harmonic function is also $P$-harmonic. Although the
converse is not true in general, one can show that \emph{any function which is
minimal in the cone of non-negative strongly harmonic functions is also minimal
in the cone of non-negative $P$-harmonic functions} (which, in particular,
implies coincidence of $P$-harmonicity and strong harmonicity for bounded
functions, see \remref{rem:strong} below). Thus, the above functions $\f_\ga$
are minimal $P$-harmonic, which is precisely what is needed for applying
\thmref{thm:amenab actions}.

We shall briefly outline the proof (which actually works for any Markov
operator $P$ from the ``Hecke algebra'' generated by a family of commuting
Markov operators $\{P_i\}$ under suitable non-degeneracy conditions), the
details to be given elsewhere. Let $\f$ be a minimal strongly harmonic
function. By passing from the operators $P_i$ and $P$ to their Doob transforms
(which commute simultaneously with the original operators) we may assume that
$\f=\1$. If $\1$ is not minimal $P$-harmonic, then the Poisson boundary of the
operator $P$ is non-trivial. Let $\psi$ be the bounded $P$-harmonic function
which corresponds to a non-constant measurable function $\wh\psi$ on the
Poisson boundary with the values 0 and 1 only. By the martingale convergence
theorem the values of $\psi$ converge to the function $\wh\psi$ along a.e.
sample path $(x_n)$ of the Markov chain associated with the operator $P$, so
that these limits are either 0 or 1. On the other hand,
$$
\psi(x_n) = P\psi(x_n) = \frac1d \bigl( P_1\psi(x_n) + \dots + P_d\psi(x_n)
\bigr) \;.
$$
Since the operators $P_i$ and $P$ commute, the functions $P_i\psi$ are also
$P$-harmonic and take values between 0 and 1. Therefore, their limits along
almost all sample paths are the same as for the function $\psi$. Since bounded
harmonic functions are determined by their limit values on the Poisson
boundary, we conclude that $P_i\psi=\psi$, i.e., $\psi$ is strongly harmonic in
contradiction to the hypothesis of minimality of the constant function $\1$ in
the cone of strongly harmonic functions.

\begin{rem} \label{rem:strong}
The functions $\f_\ga,\,\ga\in\part X$ provide a decomposition of the constant
function. Since, as we have just proved, $\f_\ga$ are minimal $P$-harmonic,
\emph{the space $\part X$ with the corresponding representing measure is the
Poisson boundary of the operator} $P$. In particular, \emph{for bounded
functions on $X$ strong harmonicity is equivalent to $P$-harmonicity} (cf.
\cite{Mantero-Zappa03} for the dimension 2 case).
\end{rem}

\subsection{Amenability at infinity}

A locally compact group $G$ is called \emph{amenable at infinity} if it admits
a topologically amenable action on a Hausdorff compact space $X$, i.e., if the
associated groupoid $G\lsd X$ is topologically amenable. For a countable group
$G$ its amenability at infinity is equivalent to topological amenability of its
action on the Stone-\v{C}ech compactification $\be G$, and, moreover, if $G$ is
finitely generated, to existence of a uniform embedding of its Cayley graph
into Hilbert space \cite{Higson-Roe00}. This notion has found important
applications in the theory of $C^*$ algebras, see \cite{Higson-Roe00},
\cite{Higson00}, \cite{Anantharaman02}, \cite{Valette02},
\cite{Chabert-Echterhoff-Oyono03} and the references therein.
\thmref{thm:geometric amenability} implies

\begin{thm}
Any closed subgroup of the group of isometries of any of the spaces listed in
\thmref{thm:geometric amenability} is amenable at infinity.
\end{thm}

It is known that any discrete subgroup of a connected Lie group is amenable at
infinity \cite[Example~5.2.2]{Anantharaman-Renault00}. As it follows from the
theorem of Adams on universal amenability of the boundary action of the group
of isometries of an exponentially bounded Gromov hyperbolic space
\cite{Adams96}, any discrete group of isometries of such a space is also
amenable at infinity in view of \cite[Theorem 3.3.7]{Anantharaman-Renault00}
(also see \cite{Germain00} for the particular case of word hyperbolic groups).
It was proved in \cite{Kaimanovich03} that this is true for general closed
groups of isometries of Gromov hyperbolic spaces as well under suitable bounded
geometry assumptions (without which amenability of the boundary action may
fail).

\begin{rem}
Disproving a conjecture from \cite{Higson-Roe00}, Gromov \cite{Gromov00} showed
that there exist finitely generated groups $G$ whose Cayley graph does not
admit a uniform embedding into Hilbert space, and which, therefore, are not
amenable at infinity. In view of \thmref{thm:martin} these groups have a
curious property: the set of minimal harmonic functions in the Martin boundary
of any random walk on $G$ (more generally, on any proper $G$-space) is never
closed.
\end{rem}

\bibliographystyle{amsalpha}

\enddocument

\bye